\documentclass{article}
\usepackage{graphicx}

\def\be{\begin{equation}}
\def\ee{\end{equation}}

\newcommand{\ff}[1]{{\bf  #1}}

\def\b{\beta}

\def\ra{\rightarrow}
\def\x{\ff{x}}

\begin{document}

\title{Applications and Analysis of Bio-Inspired Eagle Strategy for Engineering Optimization}
\author{Xin-She Yang and Mehmet Karamanoglu  \\
School of Science and Technology, Middlesex University, \\
The Burroughs, London NW4 4BT, United Kingdom. \\
\and
T. O. Ting \\
Department of Electrical and Electronic Engineering,  \\
Xi'an Jiaotong-Liverpool University, \\
Suzhou, Jiangsu Province, P. R. China. \\
\and
Yu-Xin Zhao \\
College of Automation, Harbin Engineering University, \\ 
Harbin, P. R. China.
}

\date{}

\maketitle

\begin{abstract}
All swarm-intelligence-based optimization algorithms use some stochastic components
to increase the diversity of solutions during the search process. Such randomization
is often represented in terms of random walks. However, it is not yet clear why
some randomization techniques (and thus why some algorithms) may perform better
than others for a given set of problems. In this work, we analyze these randomization
methods in the context of nature-inspired algorithms. We also use eagle strategy
to provide basic observations and relate step sizes and search efficiency using
Markov theory. Then, we apply our analysis and observations to solve four
design benchmarks, including the designs of a pressure vessel, a speed reducer, a PID controller and a heat exchanger.
Our results demonstrate that eagle strategy with L\'evy flights can perform
extremely well in reducing the overall computational efforts.

\end{abstract}
{\bf Citation details:}
X. S. Yang, M. Karamanoglu, T. O. Ting and Y. X. Zhao, 
Applications and Analysis of Bio-Inspired Eagle Strategy for Engineering Optimization，
{\it Neural Computing and Applications}, vol. 25, No. 2, pp. 411-420 (2014).


\section{Introduction}

In contemporary neural computing, an active branch of research is the
nature-inspired algorithms with diverse applications in engineering optimization.
Most of these algorithms are based on the so-called swarm intelligence,
and usually involve some form of non-deterministic,
stochastic components, which often appear in terms of random walks.
Such random walks can be surprisingly efficient when combined with
deterministic components and elitism, as this has been demonstrated in
many modern metaheuristic algorithms such as particle swarm optimization,
firefly algorithm and other algorithms \cite{Blum,Kirk,Ting,YangComp,YangCam,YangAPSO,YangBAE,YangGEWA,YangEagle}.

Recent studies in nature-inspired algorithms have shown promising results
with divers algorithms, including new algorithms such as accelerated particle swarm
optimization \cite{YangAPSO,GandomiAPSO}, bat algorithm \cite{YangBAE},
krill herd algorithm \cite{GandomiKH}, flower algorithm \cite{YangFPA}, and
other algorithms \cite{YangEagle,Ting}. A comprehensive review can be found in \cite{Fister,GandomiYang,YangMOCS}.
In all these algorithms, different degrees of randomization, exploration and exploitation
have been used so as to maintain a good degree of solution diversity in the solution
population, which helps to enhance the performance of these algorithms.
Applications of modern nature-inspired algorithms have been very diverse
with promising results \cite{Fister,Fister2,Ting2,GandomiYangBA,YangMOBA}.

In order to gain insight into the working mechanism of a stochastic algorithm,
mathematical analysis of the key characteristics of random walks is necessary.
Though there are some extensive studies of random walks with solid results
in the statistical literature, most of these results are based on rigorous
assumptions so as to obtain  theoretical results using Markov
chain models and/or Markov chain Monte Carlo methods \cite{Fishman,Gamer,Geyer,Ghate,Gilks,YangComp}.
Consequently, such results may be too theoretical, and thus have
not much practical implications for designing optimization algorithms.
In addition, it is necessary to translate any relevant theoretical
results in the right context so that they are truly useful
to the optimization communities. The current work has extended our earlier work extensively \cite{Yangetal2013}.
Therefore, the aims of this paper are two-folds: to introduce the random walks
and L\'evy flights in the proper context of metaheuristic optimization, and
to use these results in the framework of Markov theory
to analyze the iteration process of algorithms such as step sizes,
efficiency and the choice of some key parameters.

The rest of the paper is organized as follows: Section 2 briefly introduce the eagle strategy (ES).
Section 3 introduces the fundamentals of random walks
and discusses L\'evy flights, as well as their links to optimization
via Markov chain theories. Section 4 analyzes the choice of step sizes,
stopping criteria and efficiency.
Section 5 presents four case studies for engineering optimization applications.
Finally, we briefly draw the conclusions in Section 6.

\section{Eagle Strategy and Solution Diversity}

\subsection{Eagle Strategy}
Eagle strategy developed by Xin-She Yang and Suash Deb \cite{YangEagle} is a two-stage method
for optimization. It uses a combination of crude global search and intensive local
search employing different algorithms to suit different purposes. In essence,
the strategy first explores the search space globally using a L\'evy
flight random walk; if it finds a promising solution, then an intensive
local search is carried out by using a more efficient local optimizer such as hill-climbing and downhill
simplex method. Then, the two-stage process starts again with
new global exploration followed by a local search in a new region.
The main steps of this method can be represented as the pseudo code as outlined in Fig. \ref{eagle-fig-50}.

The advantage of such a combination is to use a balanced tradeoff between
global search which is often slow and a fast local search. Some tradeoff and balance
are important. Another advantage of this method is that we can use any algorithms we like
at different stages of the search or even at different stages of iterations.
This makes it easy to combine the advantages of various algorithms so as
to produce better results.

It is worth pointing out that this is a methodology or strategy, not an algorithm.
In fact, we can use different algorithms at different stages and at different
time of the iterations.  The algorithm used for the global exploration should
have enough randomness so as to explore the search space diversely and
effectively. This process is typically slow initially, and should speed up
as the system converges, or no better solutions can be found after a certain
number of iterations. On the other hand, the algorithm used for the intensive
local exploitation should be an efficient local optimizer. The idea is to reach the
local optimality as quickly as possible, with the minimal number of function evaluations.
This stage should be fast and efficient.

\begin{figure}
\begin{center}
\begin{minipage}[c]{0.7\textwidth}
\hrule \vspace{2pt}
Objective function $f(\x)$ \\
Initialization and random initial guess $\x^{t=0}$ \\
{\bf while} (stop criterion) \\
Global exploration by randomization (e.g. L\'evy flights) \\
Evaluate the objectives and find a promising solution \\
Intensive local search via an efficient local optimizer  \\
\indent \qquad {\bf if} (a better solution is found) \\
\indent \qquad \quad Update the current best \\
\indent \qquad {\bf end} \\
Update $t=t+1$ \\
{\bf end} \\
\hrule
\end{minipage}
\caption{Pseudo code of the eagle strategy. \label{eagle-fig-50}  }
\end{center}
\end{figure}

For the local optimizer in this paper, we will use the accelerated particle swarm optimization (APSO) \cite{YangAPSO}
which is a simple but efficient variant of particle swarm optimization.
The APSO essentially has one updating equation
\be \x_i^{t+1} =(1-\beta) \x_i^t + \beta \ff{g}^* + \alpha \epsilon_t, \ee
where $\ff{g}^*$ is the current best solution among all the solutions $\x_i^t$ at iteration $t$.
$\beta \in (0, 1) $ is a parameter, and $\alpha=O(1)$ is the scaling factor. Here $\epsilon_t$
is a random number drawn from a standard normal distribution $N(0,1)$. As the typical scale $L$
of a problem may vary, $\alpha$ should be linked to $L$ as $\alpha=0.1 \alpha_t L$ where
$\alpha_t$ should decrease as the iterations proceed in the following form
\be \alpha_t =\alpha_0 \gamma^t, \ee
where $\alpha_0 \in (0.5, 1)$ and $0<\gamma<1$. From our previous parametric study, we will
use $\alpha_0=1$, $\beta=0.5$, and $\gamma=0.97$ \cite{YangAPSO}.

\subsection{Exploration, Exploitation and Solution Diversity}

In almost all nature-inspired algorithms, two conflicting and yet important components
are exploration and exploitation, or diversification and intensification. The balance between these
components are very important to ensure the good performance of an algorithm \cite{Blum,Ferrante}.
In other words, a good degree of diversity should be maintained in the population of the solutions
so that exploration and exploitation can be reflected in the evolving population. If the population
is too diverse, it is good for global exploration, but it may slow down the convergence. On the other
hand, if the diversity is too low, intensive local exploitation may lead to premature convergence,
and thus may loose the opportunity of finding the global optimality. However, how to maintain good balance
is still an unsolved problem, and different algorithms may have different ways of dealing with this issue.

In most algorithms such as the particle swarm optimization, exploration and diversity can be considered
as the steps by using random numbers, while intensification and exploitation are by the use
of current global best. However, there is no direct control on how to switch between these components.
Even in the accelerated particle swarm optimization (ASPO), diversity is mainly controlled by a random-walk-like term, while
the control is indirectly carried out by an annealing-like reduction of randomness.
On the other hand, eagle strategy provides a direct control of these two stages/steps in
an iterative manner. First, solutions are sampled in a larger search space, and these solutions
often have high diversity. These solutions are then fed into the APSO for evolution so that
a converged state can be reached, and at the converged state, solution diversity is low.
Then, a new set of samples are drawn again from the larger search space for another round of
intensive APSO iteration stage. In this way, both exploration and exploitation have been used to main
a good degree of diversity in the overall population, which also allows the system to converge
periodically towards global optimality.

As we will see in the late analysis of random walks and search strategies, randomization techniques
are often used for exploration to increase the diversity of the solutions, while the selection of good solutions
and evolution of an algorithm tend to lead to convergence of the system. However, their role
is subtle. For example, in memetic algorithms, the balance is even more subtle \cite{Ferrante}.
It can be expected that the analysis in the rest of the paper can provide some insight into the
working mechanisms and subtlety of random walks and randomization techniques in maintaining the
good diversity of solutions in different algorithms.

\section{Random Walks and L\'evy Flights as a Search Strategy}

In modern stochastic optimization algorithms, especially those based on swarm
intelligence, there are often a deterministic component and a stochastic component, though traditional algorithms
such as the steepest descent method
are purely deterministic. Randomness is now an essential part of the stochastic
search algorithms.

Randomization techniques such as random walks have become an integrated part of
a search process in stochastic algorithms. However, how to achieve the effective
randomization remains an open question. The exact form of randomization
may depend on the actual algorithm of interest. One of the objective of this
paper is to analyze and discuss  the main concepts of random walks and L\'evy flights, and their role in metaheuristic optimization.

\subsection{Search via Random Walks}
In essence, a stochastic search process involves the use of a random process
to generate new solutions in the search space so that the solutions can sample
the landscape appropriately. However, the effectiveness of this sampling process
depends on the random process used and the actual way of generating new solutions/samples.

If the search space is treated as a black box (thus no knowledge or
assumption about the modality is made), a random walk is one of the most simplest ways
to carry out the search for optimality. Briefly speaking,
a random walk is a random process which consists of taking a
series of consecutive random steps \cite{YangComp,YangGEWA}.
That is, the total moves $S_N$ after $N$ steps are the sum of each
consecutive random step $X_i (i=1,...,N)$:
\be S_N=\sum_{i=1}^N X_i = X_1 + ... + X_N = \sum_{i=1}^{N-1} X_i + X_N = S_{N-1} + X_N, \ee
where $X_i$ is a random step drawn from a random distribution such as a normal distribution.
Depending on the perspective, the above relationship can also be
considered as a recursive formula. That is, the next state $S_N$
will only depend on the current existing state $S_{N-1}$
and the move or transition $X_N$ from the existing state to the next state.
In other words, the next state will depend only on the current state
and the transition probability, and it has no direct link to the states in the past.
From the Markov theory, we know that this is typically the
main property of a Markov chain, to be introduced later.

It is worth pointing out that there is no specific restriction on the step sizes.
In fact, the step size or length in a random walk can be fixed or varying.
Random walks have many applications in physics, economics, statistics,
computer sciences, environmental science and engineering.
Mathematically speaking, a random walk can be defined as
\be S_{t+1}= S_t + w_t, \ee
where $S_t$ is the current location or state at $t$, and $w_t$ is a step or
random variable drawn from a known probability distribution.

If each step or jump is carried out in the $d$-dimensional space, the random walk $S_N$
discussed earlier becomes a random walk in higher dimensions. In addition,
there is no reason why the step length should be fixed. In general,
the step size can also vary according to a known distribution.
If the step length obeys the Gaussian distribution, the random walk becomes the
Brownian motion or a diffusion process.

As the number of steps $N$ increases,
the central limit theorem implies that the random walk should approach
a Gaussian distribution. If the steps are drawn from a normal distribution
with zero mean, the mean of particle locations is obviously zero.
However, their variance will increase linearly with $t$. This is valid for
the case without any drift velocity.

In a more generalized case in the $d$-dimensional space,
the variance of Brownian random walks can be written as
\be \sigma^2(t) = |v_0|^2 t^2 + (2 d D) t, \ee
where $v_0$ is the drift velocity of the system \cite{YangGEWA,Viswan}.
Here $D=s^2/(2 \tau)$ is the effective diffusion coefficient
which is related to the step length $s$ over a short time interval $\tau$ during each jump.
For example, the well-known Brownian motion $B(t)$ obeys a Gaussian distribution with zero mean and time-dependent variance:
\be B(t) \sim N(0, \sigma^2(t)), \ee
where $\sim$ means the random variance obeys the distribution on the right-hand side; that is, samples should be drawn from the distribution.

In physics and chemistry, a diffusion process can be considered
as a series of Brownian motion, which obeys the Gaussian distribution.
Therefore, standard diffusion is often referred to as the Gaussian diffusion.
If the motion at each step is not Gaussian, then the diffusion is called non-Gaussian diffusion. On the other hand, if the step lengths are drawn from other distributions,
we have to deal with more generalized random walks. For example, a very special case
is when step lengths obey the L\'evy distribution, such a random walk is called  L\'evy flight or L\'evy walk \cite{Pav,Viswan}. It is worth pointing out that a L\'evy flight is also a Markov chain. In fact, any algorithmic path traced by the current solution plus a transition probability forms a Markov chain. This is one of the reason
why Markov chain theory can be used to analyzed stochastic algorithms such as simulated annealing and cuckoo search.

In Section 4, we will discuss random walks without drift. That is, we will set $v_0=0$.

\subsection{L\'evy Flights}

In standard Gaussian random walks, the steps are drawn from a Gaussian normal distribution
$N(0,\sigma)$, and these steps are mostly limited within $3 \sigma$. Therefore, it can be expected that very large steps ($>3 \sigma$) are extremely rarely. Sometimes, it may be necessary to generate new solutions that are far from the current state so as to avoid being trapped in a local region for a long time. In this case, random walks with varying step sizes may be desirable. For example, L\'evy flights are another class of random walks
whose step lengths are drawn from the so-called L\'evy distribution. When steps are
large, L\'evy distribution can be approximated as a simple power-law
\be L(s) \sim |s|^{-1-\b}, \ee
where $0<\b \le 2$ is an index \cite{Gutow,Mant,Nol,Pav}. However, this power-law is
just an approximation to the L\'evy distribution.

To be more accurate, L\'evy distribution should be defined in
terms of the following Fourier transform
\be F(k)=\exp[-A |k|^{\b}], \quad 0 < \b \le 2, \ee
where $A$ is a scaling parameter. In general, the inverse of this
integral is not straightforward, as no analytical form can be obtained,
except for a few special cases. One special case $\b=2$
corresponds to a Gaussian distribution, and another case $\b=1$
leads to a Cauchy distribution.

Though the inverse integral is difficult, however, one useful technique
is to approximate
\be L(s) =\frac{1}{\pi} \int_0^{\infty} \cos (k s) \exp[-A |k|^{\b}] dk, \ee
when $s$ is large. That is,
\be L(s) \ra \frac{A \; \b \; \Gamma(\b) \sin (\pi \b/2)}{ \pi |s|^{1+\b}},
\quad s \ra \infty, \ee
where the Gamma function $\Gamma(z)$ is defined as
\be \Gamma(z) =\int_0^{\infty} t^{z-1} e^{-t} dt.  \ee
Obviously, when $z=n$ is an integer, it becomes $\Gamma(n)=(n-1)!$.

Since the steps drawn from a L\'evy distribution can be occasionally very large,
random walks in terms of L\'evy flights are more efficient than standard
Brownian random walks. It can be expected that L\'evy flights are efficient
in exploring unknown, large-scale search space. There are many reasons
to explain this high efficiency, and one reason is that the variance of L\'evy flights
\be \sigma^2(t) \sim t^{3-\b},  \quad 1 \le \b \le 2, \ee
increases much faster than the linear relationship (i.e., $\sigma^2(t) \sim t$) of Brownian random walks.

It is worth pointing out that a power-law distribution is often linked to
some scale-free characteristics, and L\'evy flights can thus show
self-similarity and fractal behavior in the flight patterns.
Studies show that L\'evy flights can maximize the efficiency of the
resource search process in uncertain environments.
In fact,  L\'evy flights have been observed among the
foraging patterns of albatrosses, fruit flies, and spider monkeys. Even humans
such as the Ju/'hoansi hunter-gatherers can trace paths of L\'evy-flight patterns \cite{Romos,Reynolds,Reynolds2,Viswan}.
In addition, L\'evy flights have many applications.  Many physical phenomena such as
the diffusion of fluorenscent molecules, cooling behavior and noise could show L\'evy-flight
characteristics under right conditions.

L\'evy flights have been successfully used in optimization to enhance
the search efficiency
of nature-inspired algorithms \cite{YangMOCS,YangGEWA,YangMOFA,Romos}.
The above nonlinear variance
partly explain why. As we will see below, a good combination with other methods such as
eagle strategy can be even more efficient.

\subsection{Optimization as Interacting Markov Chains}

If we look at an algorithm from the Markovian view, an algorithm is
intrinsically related to Markov chains because an algorithm is
an iterative procedure whose aim is to generate new,
better solutions from the current solution set so that the best solution can be reached in a finite number of steps,
ideally, as fewer steps as possible. In this sense, the next solutions (i.e., states)
can depend only on the current solution (states) and the way to move (the transition)
towards the new solution (i.e., new states). Therefore, the solution paths are Markov chains.

In the very simplest case, a very good example is the so-called simulated annealing \cite{Kirk,YangComp}, which is a Markov chain generating a
piece-wise path in the search space.
Broadly speaking, swarm-intelligence-based algorithms such as particle swarm
optimization, bat algorithm and  eagle strategy can all be considered
as a system of multiple interacting Markov chains \cite{YangBAE,YangEagle}.
Now, let us discuss these concepts in detail.

Briefly speaking, a random variable $U$ is said to form a
Markov process if the transition probability, from state $U_t=S_i$
at time $t$ to another state $U_{t+1}=S_j$, depends
only on the current state $U_i$, independent of any past states before $t$.
Mathematically, we have
\be P(i,j) \equiv P(U_{t+1}=S_j\Big|U_0=S_{p}, ..., U_t=S_i)
 = P(U_{t+1}=S_j\Big|U_t=S_i), \ee
which is independent of the states before $t$.
The sequence of random variables $(U_0, U_1, ..., U_n)$
generated by a Markov process is subsequently
called a Markov chain. Obviously, a random walk is a Markov chain.

The transition probability
\be P(i,j) \equiv P( i \ra j)=P_{ij}, \ee
is also called the transition kernel of the Markov chain.
From the algorithmic point of view, different algorithms will have different
transition kernel; however, it is not known what kernels are most effective
for a given problem. In order to solve an optimization problem,
the feasible solution set can be obtained
by performing a random walk,
starting from a good initial but random guess solution. However, simple or
blind random walks are not efficient.

To be computationally efficient and effective in searching for new solutions,
effective transition kernels should allow to generate new solutions near the truly
optimal solutions as well as to increase the mobility of the random walk
so as to explore the search space more effectively.
In addition, the best solutions found so far should be kept in the population.
Ideally, the way to control the walk should be carried out
in such a way that it can move towards the optimal solutions more quickly,
rather than wandering away from the potential best
solutions. These are the challenges for most metaheuristic algorithms,
and various attempts are being made so as to design better optimization algorithms.

\section{Search Efficiency and Step Sizes}

\subsection{Step Sizes in Random Walks}

In all metaheuristic algorithms, different forms of random walks are widely
used for randomization and local search \cite{YangComp,YangBAE}. Obviously, a proper
step size is very important.

Many algorithms typically use the following generic equation:
\be \x^{t+1}=\x^t + s \; \ff{\epsilon}_t, \ee
where $\ff{\epsilon}_t$ is drawn from a standard normal distribution with
zero mean and unity standard deviation. Here, the step size $s$ is essentially
a scaling factor, controlling how far a random walker,
such as an agent or a particle in metaheursitics, can move for a fixed number of iterations.

From the above equation, it can be expected that the new solution $\x^{t+1}$ generated will be
too far away from the old solution (or more often the current best) if $s$ is too large.
Then, a move that is too far away is unlikely to be accepted as a better solution.
On the other hand, if $s$ is too small, the change is too small to be significant,
and the new solution may be too close to the existing solution.
Consequently, the diversity of the new solutions is limited, and thus
the search process is not efficient.
Therefore, an appropriate step size is important to maintain
the search process as efficient as possible.

There are extensive good theoretical results about isotropic random walks in
the literature \cite{Gutow,Mant,Nol,Pav,Reynolds,YangComp}, and one of the results
concerns the average distance $r$ traveled in the $d$-dimension space, which is
\be r^2=2 d D t. \ee
Here, $D=s^2/2\tau$ is the effective diffusion coefficient where $s$ is the step size
or distance traveled at each jump. In addition, $\tau$ is the time taken for each jump.
By re-arranging the above equation, we get
\be s^2=\frac{\tau \; r^2}{t \; d}, \ee
which can be used to estimate the typical step sizes for a give problem.
For example, for a typical scale $L$ of dimensions of interest, the local search
is typically limited in a region of $r=L/10$.  As the iterations are discrete,
$\tau=1$ can used for simplicity. Obviously, the number of iterations should not be too large;
otherwise, the computational costs are too high. Typically, the number
of generations is usually $t=100$ to $1000$ for most applications. Therefore, we have
\be s \approx \frac{r}{\sqrt{t d}}=\frac{L/10}{\sqrt{t \; d}}. \label{stepsize-equ-555} \ee

Let us try to do some estimates. For $d=1$ and $t=100$, we have $s \approx 0.01L$,
while $s \approx 0.001L$ for $d=10$ and $t=1000$. In addition, for $d=100$ and $t=1000$,
we have $s \approx L/3000$. As step sizes could differ from variable to variable, a step size
ratio $s/L$ is more generic. Therefore, we can use $s/L=0.001$ to $0.01$ for most applications.
In the rest of the paper, we will usually set the number of iterations
as 200 in the case studies in Section 5.

\subsection{Accuracy and Number of Iterations}

If an algorithm works well, the final accuracy of the obtained solution will depend
on the number of iterations. In principle, a higher number of iterations may
be more likely to obtain higher accuracy, though stagnation may occur. From the theory
of random walks, we can estimate the number of iterations needed for a given tolerance,
though such estimates are just guidelines. For example, in order to achieve the
accuracy of $\delta=10^{-5}$, the number of steps or iterations $N_{\max}$ needed by pure random walks can be estimated by
\be N_{\max} \approx  \frac{L^2}{\delta^2 d}, \ee
which is essentially an upper bound.  As an estimate for $L=10$ and $d=100$, we have
\be N_{\max} \approx \frac{10^2}{(10^{-5})^2 \times 100} \approx 10^{10},  \label{Nmax-old} \ee
which is a huge number that is not easily achievable in practice.
However, this number is still far smaller than that needed by a
uniform or brute force search method.
It is worth pointing out the above estimate is the upper limit for the
worst-case scenarios. In reality, most metaheuristic algorithms require far fewer
numbers of iterations.

Though the above estimate may be crude, it does imply
another interesting fact that the number of iterations
will not be affected much by dimensionality. In fact, higher-dimensional problems
do not necessarily increase the number of iterations. This may lead to
a rather surprising possibility that random walks may be efficient in
higher dimensions if the optimization problem is highly multimodal. This may provide some
hints for designing better algorithms by cleverly using random walks and other
randomization techniques. As we will see below, different random walks will indeed lead to different
convergence rates, and L\'evy flights are one of the best randomization techniques.

\subsection{Why L\'evy Flights and Eagle Strategy are so Efficient}

As mentioned earlier, the variance of a Gaussian random walk increases
linearly with time $t$, while the variance of L\'evy flights usually increases nonlinearly at a higher rate.
As a result, if L\'evy flights instead of Gaussian random walks are used, the above estimate becomes
\be N_{\max} \approx \Big(\frac{L^2}{\delta^2 d}\Big)^{1/(3-\beta)}.  \label{Levy-step} \ee
If we use $\beta=1.5$ together the same values of $L=10$, $d=100$ and $\delta=10^{-5}$, we have
\be N_{\max} \approx 4.6 \times 10^6. \label{Nmax-mid} \ee
It can be clearly seen that L\'evy flights can
reduce the number of iterations by about 4 orders [$O(10^4)$]
from $O(10^{10})$ in Eq.~(\ref{Nmax-old}) to $O(10^6)$ in Eq.~(\ref{Nmax-mid}).

Ideally, the step sizes in any nature-inspired algorithm should be
controlled in such a way that they can do both local and global search more efficiently.
To illustrate this point, let us split the search process into two stages
as those in the efficient Eagle Strategy (ES), developed by Xin-She Yang and Suash Deb \cite{YangEagle,Yang2012Eagle}.
The first stage uses a crude/large step, say, $\delta_1=10^{-2}$, and then
in the second stage, a finer step size $\delta_2=10^{-5}$ is used
so as to achieve the same final accuracy as discussed in the previous section.
The first stage can cover the whole region $L_1=L$,
while the second region should cover smaller, local regions of size $L_2$.
Typically, $L_2=O(L_1/1000)$. Using the above values and $L_1=L=10$ and $L_2=L/1000=0.01$, then
we have
\be N_{1,\max} \approx 10^{4}, \quad N_{2,\max} \approx 10^{4}. \label{Nmax-new} \ee
It can be seen clearly that the number of iterations can be reduced by about 6 orders ($10^{6}$) from $O(10^{10})$ [see Eq.~(\ref{Nmax-old})]
to $O(10^4)$ [see Eq.~(\ref{Nmax-new})]. This approximate values for the number of iterations
have been observed and used in the literature \cite{YangMOCS,GandomiYangBA}. For example, the typical
number of iterations for bat algorithm in \cite{GandomiYangBA} and multiobjective cuckoo search in \cite{YangMOCS}
are typically 20,000, which were indeed consistent with our estimations here.

It is worth pointing out that the above reduction is by the two-stage Eagle Strategy only
without using L\'evy flights. It can be expected that L\'evy flights can reduce the number of
iterations even further. In fact,  if L\'evy flights are used within the Eagle Strategy,
then the above estimates can be reduced to
\be N_{1,\max} \approx N_{2,\max} \approx 464, \ee
which is obtained by substituting $L_1$ (or $L_2$) and $\delta_1$ (or $\delta_2$) into
Eq.~(\ref{Levy-step}).  The relative lower number of iterations
can be both practical and realistic. Therefore, the good combination
of L\'evy flights with Eagle Strategy can reduce the number of iterations from $O(10^{10})$
to less than $O(10^3)$, which works almost like a magic. This shows that, with a combination of good algorithms, eagle strategy can significantly reduce the computational efforts and may thus increase the search efficiency dramatically.
It may be possible that a multi-stage eagle strategy can be developed to
enhance this efficiency even further.

\section{Applications}

The above analyses and observations have indicated that
a proper combination of randomization techniques such as L\'evy flights,
eagle strategy and a good optimizer such as
APSO can reduce the computational efforts dramatically and thus improve the search efficiently significantly.
Henceforth, we will use eagle strategy with APSO to solve four
nonlinear benchmarks so as to validate the above theoretical results.
The four optimization benchmarks are: design optimization of a pressure vessel,
a speed reducer, a PID controller and a heat exchanger.

The parameter settings for the algorithms used in the rest of paper have been based on
a parametric study. For the eagle strategy, we used 5 rounds of two-stage iterations,
and each round used the APSO for an intensive search. The population size is $n=20$,
and the iteration for APSO was set to $t=200$. This leads to a total of $5*20*200=20,000$ function
evaluations for each case study. In addition, the parameters in APSO have been set
to be $\alpha_0=1, \beta=0.5$, and $\gamma=0.97$.

\subsection{Pressure Vessel Design}

Pressure vessels are literally everywhere such as champagne bottles and gas tanks.
For a given volume and working pressure, the basic aim of designing
a cylindrical vessel is to minimize the total cost. Typically, the design variables are
the thickness $d_1$ of the head, the thickness $d_2$ of the body, the inner radius
$r$, and the length $L$ of the cylindrical section \cite{Cagnina,GandomiYang}.
This is a well-known test problem for optimization  and it can be written as
\be \textrm{minimize } f(\x) = 0.6224 d_1 r L + 1.7781 d_2 r^2
 + 3.1661 d_1^2 L + 19.84 d_1^2 r, \ee
subject to the following constraints:
\be
\begin{array}{lll}
 g_1(\x) = -d_1 + 0.0193 r \le 0 \\
 g_2(\x) = -d_2 + 0.00954 r \le 0 \\
 g_3(\x) = - \pi r^2 L -\frac{4 \pi}{3} r^3 + 1296000 \le 0 \\
 g_4(\x) =L -240 \le 0.
\end{array}
\ee

The simple bounds are  \be 0.0625 \le d_1, d_2 \le 99 \times 0.0625, \ee
and \be 10.0 \le r, \quad L \le 200.0. \ee
Recently, Cagnina et al \cite{Cagnina} used an efficient particle swarm optimiser
to solve this problem and they found the best solution
$f_* \approx 6059.714$
at \be \x_* \approx (0.8125, \; 0.4375, \; 42.0984, \; 176.6366). \ee
This means the lowest price is about $\$6059.71$.

Using ES with APSO, we obtained the same results, but we used significantly
fewer function evaluations, comparing with APSO alone and other methods.
In fact, we use $5$ stages and each stage with a total of $200$ iterations.
This is at least 10 times less than the iterations needed by standard PSO.
This again confirmed that ES is indeed very efficient.

\subsection{Speed Reducer Design}

Another important benchmark is the design of a speed reducer
which is commonly used in many mechanisms such as
a gearbox \cite{Cagnina,GandomiYang}.
This problem involves the optimization of 7 variables, including
the face width, the number of teeth,
the diameter of the shaft and others. All variables are continuous within some limits,
except $x_3$ which only takes integer values.

\[ f(\x) = 0.7854 x_1 x_2^2 (3.3333 x_3^2+14.9334 x_3-43.0934) \]
\be -1.508 x_1 (x_6^2+x_7^2)+7.4777 (x_6^3+x_7^3)   +0.7854 (x_4 x_6^2+x_5 x_7^2) \ee

\be g_1(\x) = \frac{27}{x_1 x_2^2 x_3}-1 \le 0, \quad
g_2(\x) = \frac{397.5}{x_1 x_2^2 x_3^2}-1 \le 0, \ee

\be g_3(\x)=\frac{1.93 x_4^3}{x_2 x_3 x_6^4} - 1 \le 0, \quad
g_4(\x)=\frac{1.93 x_5^3}{x_2 x_3 x_7^4} - 1 \le 0, \ee

\be g_5(\x) =\frac{1.0}{110 x_6^3} \sqrt{\Big(\frac{745.0 x_4}{x_2 x_3} \Big)^2+16.9 \times 10^6} -1 \le 0, \ee

\be g_6(\x) =\frac{1.0}{85 x_7^3} \sqrt{\Big(\frac{745.0 x_5}{x_2 x_3} \Big)^2+157.5 \times 10^6} -1 \le 0, \ee
\be g_7(\x) = \frac{x_2 x_3}{40} -1 \le 0, \quad
 g_8(\x) = \frac{5 x_2}{x_1} -1 \le 0, \ee
\be g_9(\x)=\frac{x_1}{12 x_2} -1 \le 0, \quad
 g_{10}(\x) =\frac{1.5 x_6+1.9}{x_4}-1 \le 0, \ee
\be g_{11}(\x) =\frac{1.1 x_7+1.9}{x_5}-1 \le 0, \ee
where the simple bounds are $ 2.6 \le x_1 \le 3.6$, $0.7 \le x_2 \le 0.8$,
$17 \le x_3 \le 28$, $7.3 \le x_4 \le 8.3$, $7.8 \le x_5 \le 8.4$, $2.9 \le x_6 \le 3.9$,
and $5.0 \le x_7 \le 5.5$.
In one of latest studies, Cagnina et al. \cite{Cagnina} obtained
the following solution
\be \x_*=(3.5, 0.7, 17, 7.3, 7.8, 3.350214,5.286683) \ee
with $f_{\min}=2996.348165.$

Using our ES with APSO, we have obtained the new best solution
\be \x_*=(3.5, 0.7, 17, 7.3, 7.8, 3.34336449,5.285351) \ee
with the best objective
$f_{\min}=2993.7495888$.
In existing studies \cite{Cagnina}, the total number of function evaluations
was often 350,000. Here, we used $n=20$ for $200$ iterations and 5 rounds of
eagle strategy runs, giving a total of $20,000$ function evaluations.
This saves about  93\% of the computational costs. We can see that ES not only provides better solutions but also finds
solutions more efficiently using fewer function evaluations. Again the number of
iterations (20,000) is consistent with our theoretical estimations given earlier in Section 4.

\subsection{Design of a PID Controller}
Let us use ES with APSO \cite{YangAPSO} to
design a well-known PID controller \cite{Xue,Matlab}
\be u(t)=K_p \Big[e(t) +\frac{1}{T_i} \int_0^t e(\tau) d\tau + T_d \frac{d e(t)}{dt} \Big], \ee
where $e(t)=r(t)-y(t)$ is the error signal between the response $y(t)$ and the reference input $r(t)$, and
$u(t)$ is the input signal to the so-called plant model. The well-established
Ziegler-Nichols tuning scheme can usually produce very good results. Here, we use ES with APSO
to minimize the rise time, the overshoot and settling time.

It is required to tune a third-order system with the following transfer function
\be G(s) = \frac{7}{s^3+3 s^2 + 3 s +1}, \ee
so that responses of the closed loop system to track a reference meet the
following requirements: rise time is less than 1.5 s, settling time is less than
5.5 s, and the overshoot is less than 5\%.
If we use the standard Matlab control system toolbox \cite{Matlab}, the obtained design by the Ziegler-Nichols scheme is
\be G_{PID}(s)=0.311 (1+ \frac{1}{2.3643 s} + 0.5911 s), \ee
which gives a rise time of 1.6 s, a settling time of 4.95, and overshoot
of about 8.54\%. Clearly, not all two requirements are met.

\begin{figure}
\centerline{\includegraphics[height=2in,width=3in]{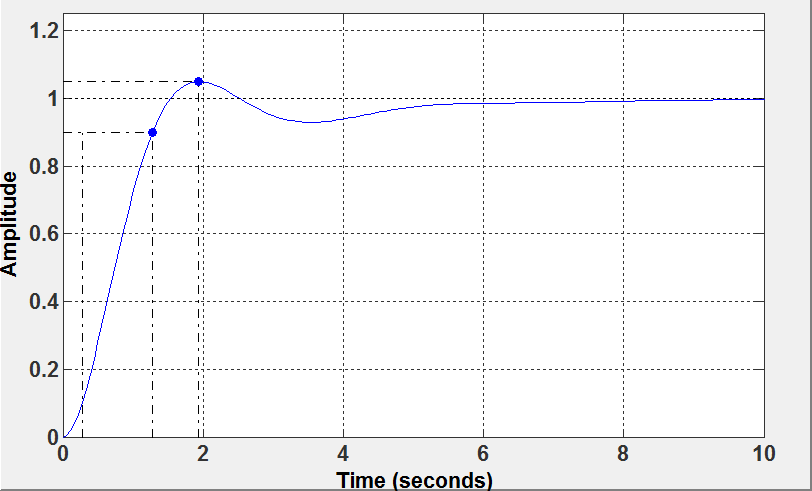}}
\caption{PID controller design by ES with APSO. \label{fig-100}}
\end{figure}

In order to refine the designs, we use ES with APSO, and the results are shown in
Fig. \ref{fig-100}. The final design becomes
\be G_{PID,new}(s)=0.5366 \Big(1 + \frac{1}{3.3940 s}+0.8485 s), \ee
which gives a rise time of 1.0 s, a settling time of 5.25 s, and
the overshoot of 4.97\%. All the design requirements are met by this new design.
As before,  the total number of function evaluations is 20,000.

\subsection{Heat Exchanger Design}
Heat transfer management is very important in many applications such as
central heating systems and microelectronics. There are many well-known benchmarks
for testing design/optimization tools, and one of the benchmarks is the design of
a heat exchanger \cite{YangBAE}.  This design problem has 8 design variables and 6 inequality constraints,
which can be written as
\be \textrm{Minimize } \; f(\x) =x_1 +x_2 +x_3, \ee
subject to
\be g_1(\x)=0.0025(x_4+x_6)-1 \le 0, \ee
\be g_2(\x)=0.0025(x_5+x_7-x_4) -1 \le 0, \ee
\be g_3(\x)=0.01(x_8-x_5) -1 \le 0, \ee
\be g_4(\x)=833.33252 x_4 +100 x_1 -x_1 x_6 - 83333.333 \le 0, \ee
\be g_5(\x)=1250 x_5 + x_2 (x_4 - x_7) -120 x_4 \le 0, \ee
\be g_6(\x)=x_3 x_5 -2500 x_5 - x_3 x_8 + 1.25 \times 10^6 \le 0. \ee
The first three constraints are linear, while the last three constraints are nonlinear.

By using ES with APSO for the same parameter values of $n=20$, $t=200$ iterations
and $5$ ES stages, we obtained the best solution
\[ \x_*=(579.30675, 1359.97076, 5109.97052, 182.01770, \]
\be \qquad 295.60118, 217.98230, 286.41653, 395.60118), \ee
which gives the optimal objective of $f_{\min}=7049.248$. This is exactly the same
as the best solution found  by Yang and Gandomi \cite{YangBAE} and is better
than the best solutions reported in the previous literature \cite{Jaber,Deb}.
In the study by Jaberipour and Khorram, they used 200,000 function evaluations,
while Deb used 320,080 function evaluations\cite{Deb}.
In the present study, we have used 20,000 function evaluations that is less than 10\% of
the computational costs by other researchers.

As we can see from the above 4 case studies, ES with APSO can save about 90\% of the
computational costs, which demonstrates the superior performance of ES with APSO.
In order to show that the improvements are significant, we use the standard
Student $t$-test in terms of the numbers of functional evaluations. For $\alpha=0.05$, the two-sample
$t$-test gives $p=0.000207$, which means that the improvements are statistically significant.

\section{Conclusions}

All swarm-intelligence-based algorithms such as PSO and firefly algorithm
can be viewed in a unified framework of Markov chains; however, theoretical analysis
remains challenging. We have used the fundamental concepts of random walks and
L\'evy flights to analyze the efficiency of random walks in
nature-inspired metaheuristic algorithms.

We have demonstrated that L\'evy flights can be significantly more
efficient than standard Gaussian random walks under appropriate conditions. By the right combination with
Eagle Strategy, significant computational efforts can be saved, as we have
shown in the paper. The theory of interacting Markov chains is complicated
and yet still under active development; however, any progress in such areas will play a central role in understanding how population- and trajectory-based metaheuristic algorithms perform under various conditions.

Even though we do not fully understand why metaheuristic algorithms work, this
does not hinder us to use these algorithms efficiently. On the contrary, such
mysteries can drive and motivate us to pursue further research and development in metaheuristics. Further research can focus on the extensive testing of metaheuristics
over a wide range of large-scale problems. In addition, various statistical measures
and self-adjusting random walks can be used to improve the efficiency of existing
metaheuristic algorithms.

On the other hand, the present results are mainly concerned with Gaussian random walks, L\'evy flights and eagle strategy. It can be expected these results may be further improved with parameter tuning and parameter control in metaheuristic algorithms. It is a known fact
that the settings of algorithm-dependent parameters can influence the convergence
behaviour of a given algorithm, but how to find the optimal setting remains an open question.
It can be very useful to carry out more research in this important area.

\end{document}